\documentclass[12pt]{amsart}
\usepackage{latexsym, amssymb, amsmath}

\newtheorem{conj}{Conjecture}
\newtheorem{thm}{Theorem}[section]
\newtheorem{lem}[thm]{Lemma}
\newtheorem{cor}[thm]{Corollary}
\newtheorem{prop}[thm]{Proposition}
\newtheorem{rem}[thm]{Remark}
\theoremstyle{definition}
\newtheorem{defn}[thm]{Definition}

\newtheorem*{ex}{Example}

\newcommand{\onabla}{\overline\nabla}

\newcommand{\lapla}{\Delta}

\newcommand{\p}{\phi}

\newcommand{\met}{\langle \cdot , \cdot \rangle}

\title{Biharmonic submanifolds in manifolds with bounded curvature}

\author{Shun Maeta} 

\thanks{Supported by the Grant-in-Aid for Research Activity Start-up, No. 25887044, 
  Japan Society for the Promotion of Science. }

\address{\footnotesize{Division of Mathematics, Shimane University, Nishikawatsu 1060 Matsue, 690-8504, Japan. }
 }
\footnotesize{
\email{shun.maeta@gmail.com~{\it or}~maeta@riko.shimane-u.ac.jp}
}

\begin{document} 
\maketitle 
\markboth{Biharmonic submanifolds in manifolds with bounded curvature} 
{Shun Maeta}

\begin{abstract} 
We consider a complete biharmonic submanifold $\phi:(M,g)\rightarrow (N,h)$ in a Riemannian manifold with sectional curvature bounded from above by a non-negative constant $c$.
Assume that the mean curvature is bounded from below by $\sqrt c$. 
If 
(i) $\int_M (|{\bf H}|^2-c)^{ p }dv_g<\infty$, for some $0<p<\infty$,
or
(ii) the Ricci curvature of $M$ is bounded from below, 
then the mean curvature is $\sqrt c$.
Furthermore, if $M$ is compact, then we obtain the same result without the assumption (i) or (ii).
\end{abstract}


\qquad\\


\section{Introduction}\label{intro} 

In 1983, J. Eells and L. Lemaire \cite{jell1} proposed the problem to consider biharmonic maps.
Biharmonic maps are, by definition, a generalization of harmonic maps.
As well known, harmonic maps have been applied into various fields in differential geometry.
In 1964, J. Eells and J. H. Sampson considered the existence problem of harmonic maps between compact Riemannian manifolds.
 They showed that any continuous map from a compact Riemannian manifold into a compact Riemannian manifold of non-positive curvature is free homotopically deformable to harmonic maps. 
 By using the existence and properties of harmonic maps, one can study the structure of Riemannian manifolds.
 On the other hand, non-existence results for harmonic maps are also known.
For example, a map of degree $\pm 1$ from a $2$-dimensional torus into a $2$-dimensional sphere is not homotopic to any harmonic map.
Therefore a generalization of harmonic maps seems an important subject instead.
So far, it seems a biharmonic map.
We would like to show the existence theorem of biharmonic maps into a Riemannian manifold with positive curvature. 
Therefore, we first should consider biharmonic maps into a sphere.

G. Y. Jiang \cite{jg1} considered a biharmonic isometric immersion $\p:(M,g)\rightarrow (N,h)$
from an $m$-dimensional Riemannian manifold $(M,g)$ into an $n$-dimensional Riemannian manifold $(N,h)$.
Here, we called $\p:(M,g)\rightarrow (N,h)$ is a biharmonic submanifold, if $\p$ is a biharmonic isometric immersion.
In \cite{jg1}, he also gave some examples of non-minimal (non-harmonic) biharmonic submanifolds in $S^n$ as follows.

\vspace{10pt}

\begin{ex}
(i) and (ii) are non-minimal biharmonic submanifolds in $S^n(1)$:\\
(i) $S^{n-1}(\frac{1}{\sqrt2})\subset S^n(1),$ \\
 (ii) $S^{n-p}(\frac{1}{\sqrt2})\times S^{p-1}(\frac{1}{\sqrt2})\subset S^n(1)$, with $n-p\not=p-1$. 
\end{ex}

\vspace{10pt}

After that there are many studies of biharmonic submanifolds in spheres
(cf. \cite{16}$\sim$\cite{30},\ \cite{48},\ \cite{66},\ \cite{77},\ \cite{78},\ \cite{90},\ \cite{elco2},\ \cite{133}, etc...).
Interestingly, their examples and classification results suggest that ``any biharmonic submanifold in spheres has constant mean curvature". 
With these understandings, Balmus, Montaldo and Oniciuc \cite{21} raised the following problem.

\vspace{10pt}

\begin{conj}
Any biharmonic submanifold in spheres has constant mean curvature.
\end{conj}

\vspace{10pt}

In this paper, we call this conjecture {\em BMO conjecture}.
There are affirmative partial answers to BMO conjecture, if $M$ is one of the following:

(i) A compact biharmonic hypersurface with nowhere zero mean curvature vector field and $|B|^2\geq m$ or $|B|^2\leq m$, where $|B|^2$ is the squared norm of the second fundamental form (cf. \cite{48},\  \cite{16}).
 
(ii) An orientable biharmonic Dupin hypersurface (cf. \cite{16}).

Balmus and Oniciuc also showed that (cf. \cite{23}) :
{\it
Let $M$ be a compact non-minimal biharmonic submanifold of $S^n$. Then either 
(i) there exists a point $p\in M$ such that $|{\bf H}(p)|<1$, 
(ii) $|{\bf H}|=1.$ In this case, $M$ is a minimal submanifold of a small hypersphere $S^{n-1}(\frac{1}{\sqrt{2}})\subset S^n$.}

On the other hand, 
since there is no assumption of {\it completeness} for submanifolds in BMO conjecture, 
in a sense it is a problem in {\it local} differential geometry.  
In this paper, we reformulate BMO conjecture into a problem 
in {\it global} differential geometry as the following:  

\begin{conj} 
Any {\rm complete} biharmonic submanifold in spheres has constant mean curvature.
\end{conj} 

In this paper, we give affirmative partial answers to BMO conjecture.
Furthermore, as mentioned above, since we would like to consider biharmonic maps into a Riemannian manifold with positive curvature, we consider biharmonic isometric immersions into a Riemannian manifold with bounded curvature by non-negative constant.

The remaining sections are organized as follows. 
Section~$\ref{Pre}$ contains some necessary definitions and preliminary geometric results.
 In section~$\ref{biminimal}$, we recall biminimal submanifolds. We also show that any compact biharmonic submanifold with the mean curvature is bounded from below by $\sqrt{c}$ in a Riemannian manifold with bounded curvature from above by $c$ has constant mean curvature $\sqrt{c}$.
  We also show that any complete biharmonic submanifold with Ricci curvature bounded from below and the mean curvature is bounded from below by $\sqrt{c}$ in a Riemannian manifold with bounded curvature from above by $c$ has constant mean curvature $\sqrt{c}$.
 In section~$\ref{main}$, we show that any complete biharmonic submanifold $M$ with the mean curvature is bounded from below by $\sqrt{c}$ and $\int_M (|{\bf H}|^2-c)^{ p }dv_g<\infty$, for some $0<p<\infty$ in a Riemannian manifold with bounded curvature from above by $c$ has the constant mean curvature $\sqrt{c}$. 
In section~$\ref{other}$, we give other affirmative partial answers to BMO conjecture.
In section~$\ref{sec bicon}$, we consider a biconservative hypersurface in a space form.

\qquad\\

\section{Preliminaries}\label{Pre} 

In this section, we shall give the definitions of harmonic maps and biharmonic maps.
We also recall biharmonic submanifolds.

Let $(M,g)$ be an $m$-dimensional Riemannian manifold and $(N,h)$, an $n$-dimensional Riemannian manifold, respectively.
We denote by $\nabla$ and $\nabla^N$, the Levi-Civita connections on $(M,g)$ and $(N,h)$, respectively and by $\onabla$ the induced connection on $\p^{-1}TN$.

\vspace{10pt}

Let us recall the definition of a harmonic map $\p:(M,g)\rightarrow (N,h)$.
For a smooth map $\phi:(M,g)\rightarrow (N,h)$, the {\em energy} of $\phi$ is defined by
$$E(\phi) =\frac{1}{2}\int_M|d\phi|^2 dv_g,$$
where $dv_g$ is the volume element of $g$.
The Euler-Lagrange equation of $E$ is 
$$\tau(\p)=\displaystyle \sum^m_{i=1}\{\onabla_{e_i}d\p(e_i)-d\p(\nabla_{e_i}e_i)\}=0,$$
where $\tau(\p)$ is called the {\em tension field} of $\p$ and $\{e_i\}_{i=1}^m$ is an orthonormal frame field on $M$.
 A map $\p:(M,g)\rightarrow (N,h)$ is called a {\em harmonic map} if $\tau(\p)=0$. 

\vspace{10pt}

In 1983, J. Eells and L. Lemaire \cite{jell1} proposed the problem to consider biharmonic maps which are critical points of the bi-energy functional on the space of smooth maps between two Riemannian manifolds.
In 1986, G. Y. Jiang \cite{jg1} derived the first and the second variational formulas of the bi-energy and studied biharmonic maps.
For a smooth map $\phi:(M,g)\rightarrow (N,h)$, the {\em bi-energy} of $\phi$ is defined by
$$E_2 (\phi )=\frac{1}{2}\int_M |\tau (\phi)| ^2 dv_g.$$
The Euler-Lagrange equation of $E_2$ is 
\begin{equation}\label{NSbi}
\tau_2(\phi)=-\Delta^{\phi} \tau (\phi ) -\sum^m_{i=1} R^N (\tau (\phi )  , d\phi (e_i))d\phi (e_i)=0,
\end{equation}
 where $\Delta^{\phi}:=\displaystyle\sum^m_{i=1}\left(\onabla_{e_i}\onabla_{e_i}-\onabla_{\nabla_{e_i}e_i}\right)$ and $R^N$ is the Riemannian curvature tensor of $(N,h)$ given by $R^N(X,Y)Z=\nabla^N_X\nabla^N_YZ-\nabla^N_Y\nabla^N_XZ-\nabla^N_{[X,Y]}Z$ for $X,\ Y,\ Z\in \frak{X}(N)$.
$\tau_2(\p)$ is called the {\em bi-tension field} of $\p$.
 A map $\p:(M,g)\rightarrow (N,h)$ is called a {\em biharmonic map} if $\tau_2(\p)=0$. 

\vspace{10pt}

We also recall biharmonic submanifolds.

Let $(N,h)$ be an $n$-dimensional Riemannian manifold. Let $\p:(M^m,g)\rightarrow (N^n,h=\met)$ be an isometric immersion from an $m$-dimensional Riemannian manifold with induced metric $g=\p^{-1}h.$
In this case, we identify $d\p(X)$ with $X\in \frak{X}(M)$ for each $x\in M.$
The Gauss and Weingarten formulas are given by
\begin{equation}
\nabla^N_XY=\nabla _XY+B(X,Y),\ \ \ \ X,Y\in \frak{X}(M),
\end{equation}
\begin{equation}\label{2.Wformula}
\nabla^N_X \xi =-A_{\xi}X+\nabla^{\perp}_X{\xi},\ \ \ X\in \frak{X}(M),\  \xi \in \frak{X}(M)^{\perp},  
\end{equation}
where $B$ is the second fundamental form of $M$ in $N$, $A_{\xi}$ is the shape operator for a unit normal vector field $\xi$ on $M,$ and $\nabla^{\perp}$ denotes the normal connection on the normal bundle of $M$ in $N$.
It is well known that $B$ and $A$ are related by
\begin{equation}\label{2.BA rel}
\langle B(X,Y), \xi \rangle=\langle A_{\xi}X,Y \rangle.
\end{equation}

For any $x \in M$, let $\{e_1, \cdots, e_n\}$ be an orthonormal basis of $N$ at $x$ such that $\{e_1, \cdots, e_m\}$ is an orthonormal basis of $T_xM$. 
The mean curvature vector field ${\bf H}$ of $M$ at $x$ is also given by 
$$ 
{\bf H}(x) = \frac{1}{m} \sum_{i = 1}^m B(e_i, e_i).
$$ 

If an isometric immersion $\p:(M,g)\rightarrow (N,h)$ is biharmonic, then $M$ is called a {\em biharmonic submanifold} in $N$.
 In this case, we remark that the tension field $\tau(\p)$ of $\p$ is written as $\tau(\p)=m{\bf H}$, where ${\bf H}$ is the mean curvature vector field of $M$.
The necessary and sufficient condition for $M$ in $N$ to be biharmonic is the following: 
\begin{equation}\label{NS bih sub}
\Delta^{\p}{\bf H}+\sum_{i=1}^mR^N({\bf H},d\p(e_i))d\p(e_i)=0.
\end{equation}

From $(\ref{NS bih sub})$, by an elementally argument, we obtain the necessary and sufficient condition for $M$ in $N$ 
to be biharmonic as follows (cf. \cite{smhu1}): 
\begin{align}\label{NS bih sub separate}
\ \ \Delta^{\perp} {\bf H} - \sum_{i=1}^m B(A_{\bf H}e_i, e_i) +\left[\sum^m_{i=1} R^N( {\bf H} , d\p(e_i))d\p(e_i)\right]^{\perp} = 0, 
\end{align}
\begin{align}\label{NS bih sub separate 2}
{\rm Trace}_g\left(\nabla A_{\bf H}\right)+{\rm Trace}_g\left(A_{\nabla^{\perp}_{\cdot}{\bf H}}(\cdot)\right)-\left[\sum_{i=1}^mR^N({\bf H},e_i)e_i\right]^{T}=0,
\end{align} 
where $\Delta^{\perp}$ is the (non-positive) Laplace operator associated with the normal connection $\nabla^{\perp}$.

\begin{rem}
Biharmonic submanifolds satisfy an overdetermined problem.  $($see also \cite{KU2014}$)$.
\end{rem}

\qquad\\

\section{Biminimal submanifolds}\label{biminimal}
In this section, we recall biminimal submanifolds and show that any compact biharmonic submanifold with the mean curvature is bounded from below by $1$ in a sphere has constant mean curvature $1$. 
We also show that any complete biharmonic submanifold with Ricci curvature bounded from below and the mean curvature is bounded from below by $1$ in a sphere has constant mean curvature $1$.

E. Loubeau and S. Montaldo \cite{elsm1} introduced the notion of {\em biminimal immersion} as follows:

\begin{defn}[\cite{elsm1}]
An immersion $\p:(M^m,g)\rightarrow (N^n,h)$, $m\leq n$ is called {\em biminimal} if it is a critical point of the functional 
$$E_{2,\lambda}(\p)=E_2(\p)+\lambda E(\p),\ \ \lambda\in \mathbb{R},$$
 for any smooth variation $\{\p_{t}\}$ of the map $\p$, $\p_0=\p$ such that $\left. V=\frac{d\p_t}{dt}\right |_{t=0}$ is normal to $\p(M)$.
\end{defn}

The Euler-Lagrange equation of $E_{2,\lambda}$ is
$$[\tau_2(\p)]^{\perp}+\lambda[\tau(\p)]^{\perp}=0,$$
where $[\cdot]^{\perp}$ denotes the normal component of $[\cdot]$.
We call an immersion {\em free biminimal} if it satisfies the biminimal condition for $\lambda=0$.
 If $\p:(M,g)\rightarrow (N,h)$ is an isometric immersion, then the biminimal condition is
\begin{align}\label{biminimal eq}
\left[-\lapla^{\p}{\bf H}-\sum^m_{i=1}R^{N}({\bf H},d\p(e_i))d\p(e_i)\right]^{\perp}+\lambda{\bf H}=0,
\end{align}
for some $\lambda \in\mathbb{R}$, and then $M$ is called a {\em biminimal submanifold} in $N$.
 If $M$ is a biminimal submanifold with $\lambda \geq0$ in $N$, then $M$ is called a {\em non-negative biminimal submanifold} in $N$.

\vspace{10pt}

\begin{rem}\label{biharmonic is biminimal}
We remark that {\bf every biharmonic submanifold is free biminimal}.
\end{rem}

\vspace{10pt}

From $(\ref{NS bih sub separate})$ and  $(\ref{biminimal eq})$, we obtain the necessary and sufficient condition for $M$ in $N$ to be biminimal as follows:
\begin{align}\label{N-S biminimal}
\ \ \Delta^{\perp} {\bf H} - \sum_{i=1}^m B(A_{\bf H}e_i, e_i) +\left[\sum^m_{i=1} R^N( {\bf H} , d\p(e_i))d\p(e_i)\right]^{\perp} = \lambda{\bf H}.
\end{align} 

By using $(\ref{N-S biminimal})$, we have the following lemma.

\vspace{10pt}

\begin{lem}\label{key lemma}
Let $(N,h)$ be a Riemannian manifold with sectional curvature bounded from above by a non-negative constant $c$.
Let $\p:(M,g)\rightarrow (N,h)$ be a biminimal submanifold with $\lambda\leq2mc$ in $N$.
If $|{\bf H}|\geq \sqrt {c-\frac{\lambda}{2m}}$,
then the following inequality folds:
$$\Delta(|{\bf H}|^2-\overline{C})\geq 2|\nabla^{\perp}{\bf H}|^2
+2m (|{\bf H}|^2-\overline{C})^2,$$
where $\overline{C}=c-\frac{\lambda}{2m}.$
\end{lem}

\vspace{10pt}

\begin{proof}
By using $(\ref{N-S biminimal})$, we have
\begin{equation}
\begin{aligned}
\Delta(|{\bf H}|^2-\overline{C})
=&2|\nabla^{\perp}{\bf H}|^2+2\langle B(A_{{\bf H}}e_i,e_i),{\bf H}\rangle
-2\langle R^N({\bf H},e_i)e_i, {\bf H}\rangle
+ \lambda|{\bf H}|^2\\
=&2|\nabla^{\perp}{\bf H}|^2+2 |A_{{\bf H}}|^2
-2\langle R^N({\bf H},e_i)e_i, {\bf H} \rangle 
+ \lambda|{\bf H}|^2\\ 
\geq&
2|\nabla^{\perp}{\bf H}|^2+2m|{\bf H}|^4
-2mc|{\bf H}|^2
+ \lambda|{\bf H}|^2\\ 
\geq&
2|\nabla^{\perp}{\bf H}|^2
+2m (|{\bf H}|^2-\overline{C})^2,\\ 
\end{aligned}
\end{equation}
where the first inequality follows from the sectional curvature of $N$ is bounded from above by a non-negative constant $c$ and $|A_{\bf H}|^2\geq m|{\bf H}|^4$,
since if one considers the tensor $\Phi=A_{\bf H}-|{\bf H}|^2 Id$,
 then $|\Phi|^2=|A_{\bf H}|^2-m|{\bf H}|^4\geq0.$
\end{proof}

\vspace{10pt}

\begin{rem}
 If $M^m$ is a hypersurface in $N^{m+1}$, then
$$\langle R^N({\bf H},e_i)e_i , {\bf H} \rangle=|{\bf H}|^2 {\rm Ric}^N(\xi,\xi)$$
where ${\rm Ric}^N$ is the Ricci curvature of $N$ and $\xi\in TM^{\perp}$ is a unit normal vector field.
From this, the assumption that the sectional curvature is bounded from above by a non-negative constant $c$ can be changed as the Ricci curvature is bounded from above by a non-negative constant $c$.
If $M^m$ is a hypersurface in $N^{m+1}$, then after-mentioned theorems can also be changed the assumption that the sectional curvature is bounded from above by a non-negative constant $c$ to the Ricci curvature is bounded from above by a non-negative constant $c$.
\end{rem}

\vspace{10pt}

If $M$ is compact, applying the standard maximum principle to the elliptic inequality $\Delta(|{\bf H}|^2-\overline{C})\geq 2m (|{\bf H}|^2-\overline{C})^2$ in Lemma $\ref{key lemma}$, we obtain the following theorem.

\vspace{10pt}

\begin{thm}\label{cpt case}
Let $(N,h)$ be a Riemannian manifold with sectional curvature bounded from above by a non-negative constant $c$.
Let $\p:(M,g)\rightarrow (N,h)$ be a compact biminimal submanifold with $\lambda\leq2mc$ in $N$.
 If $|{\bf H}|\geq \sqrt {c-\frac{\lambda}{2m}}$,
  then the mean curvature is $\sqrt {c-\frac{\lambda}{2m}}$.
\end{thm}

\vspace{10pt}

\begin{proof}
By Lemma $\ref{key lemma}$, we have
$$\Delta(|{\bf H}|^2-\overline{C})\geq 2|\nabla^{\perp}{\bf H}|^2
+2m (|{\bf H}|^2-\overline{C})^2 \geq 2m (|{\bf H}|^2-\overline{C})^2.$$
Since $|{\bf H}|^2-\overline{C}\geq 0,$ by using the standard maximum principle,
 we obtain $|{\bf H}|^2-\overline{C}=|{\bf H}|^2-\left(c-\frac{\lambda}{2m}\right)=0$, that is, the mean curvature is $\sqrt {c-\frac{\lambda}{2m}}$.
\end{proof}

\vspace{10pt}

By Theorem $\ref{cpt case}$, we obtain the following corollary.

\vspace{10pt}

\begin{cor}
Let $\p:(M,g)\rightarrow (N,h)$ be a compact biharmonic submanifold in a Riemannian manifold $N$ with sectional curvature bounded from above by a non-negative constant $c$.
 If $|{\bf H}|\geq \sqrt {c}$,
  then the mean curvature is $\sqrt {c}$.
\end{cor}

\vspace{10pt}

We recall Omori-Yau's generalized maximum principle.

\vspace{10pt}

\begin{lem}[\cite{Cheng-Yau}]\label{Cheng-Yau}
Let $(M, g)$ be a complete Riemannian manifold
whose Ricci curvature is bounded from below. 
Let $u$ be a smooth non-negative function on $M$. 
Assume that there exists a positive constant $k > 0$ such that
\begin{equation}\label{GMP}
\Delta u \geq k u^2\quad {\rm on}\ \ M.
\end{equation} 
Then, $u = 0$ on $M$.
\end{lem} 

\vspace{10pt}

By using Lemma $\ref{Cheng-Yau}$ and Lemma $\ref{key lemma}$, we obtain the following theorem.

\vspace{10pt}

\begin{thm}\label{Ric below}
Let $(N,h)$ be a Riemannian manifold with sectional curvature bounded from above by a non-negative constant $c$.
Let $\p:(M,g)\rightarrow (N,h)$ be a complete biminimal submanifold with $\lambda\leq2mc$ and Ricci curvature bounded from below in $N$.
If $|{\bf H}|\geq \sqrt {c-\frac{\lambda}{2m}}$,
  then the mean curvature is $\sqrt {c-\frac{\lambda}{2m}}$.
 \end{thm}

\vspace{10pt}

By Theorem $\ref{Ric below}$, we obtain the following corollary.
 
 \begin{cor}
Let $\p:(M,g)\rightarrow (N,h)$ be a complete biharmonic submanifold with Ricci curvature bounded from below in a Riemannian manifold with sectional curvature bounded from above by a non-negative constant $c$.
If $|{\bf H}|\geq \sqrt {c}$,
  then the mean curvature is $\sqrt {c}$.
\end{cor}

\vspace{10pt}

These are affirmative partial answers to BMO conjecture.

\qquad\\


\section{Biminimal submanifolds with finite condition}\label{main}

In this section, we show the following theorem.

\vspace{10pt}

\begin{thm}\label{finite condition}
Let $(N,h)$ be a Riemannian manifold with sectional curvature bounded from above by a non-negative constant $c$.
Let $\p:(M,g)\rightarrow (N,h)$ be a complete biminimal submanifold with $\lambda\leq2mc$ in $N$.
If $|{\bf H}|\geq \sqrt {c-\frac{\lambda}{2m}}$, and
$$\int_M \left(|{\bf H}|^2-\left(c-\frac{\lambda}{2m}\right)\right)^p<\infty,$$
for some $0<p<\infty$, 
 then the mean curvature is $\sqrt {c-\frac{\lambda}{2m}}$.
\end{thm}

\vspace{10pt}

By Theorem $\ref{finite condition}$, we obtain the following corollary.

\begin{cor}
Let $\p:(M,g)\rightarrow (N,h)$ be a complete biharmonic submanifold in a Riemannian manifold $N$ with sectional curvature bounded from above by a non-negative constant $c$.
If $|{\bf H}|\geq \sqrt {c}$, and
$$\int_M \left(|{\bf H}|^2- c\right)^p<\infty,$$
for some $0<p<\infty$, 
 then the mean curvature is $\sqrt {c}$.
 \end{cor}

\vspace{10pt}

\begin{proof}
Set $f=|{\bf H}|^2-\left(c-\frac{\lambda}{2m}\right)$. By the assumption $|{\bf H}|\geq \sqrt {c-\frac{\lambda}{2m}}$, $f$ is non-negative. 
For a fixed point $x_0\in M$, and for every $0<r<\infty,$
 we first take a cut off function $\lambda$ on $M$ satisfying that 
 \begin{equation}
\left\{
 \begin{aligned}
&0\leq\lambda(x)\leq1\ \ \ (x\in M),\\
&\lambda(x)=1\ \ \ \ \ \ \ \ \ (x\in B_r(x_0)),\\
&\lambda(x)=0\ \ \ \ \ \ \ \ \ (x\not\in B_{2r}(x_0)),\\
&|\nabla\lambda|\leq\frac{C}{r}\ \ \ \ \ \ \ (x\in M),\ \ \ \text{for some constant $C$ independent of $r$},
\end{aligned} 
\right.
\end{equation}
where $B_r(x_0)$ and  $B_{2r}(x_0)$ are the balls centered at a fixed point $x_0\in M$ with radius $r$ and $2r$ respectively (cf. \cite{N-U-2},\ \cite{M}).

Let $a$ be a positive constant to be determined later. Let $b=\frac{(a+2)(a+3-d)}{a+2-d}$, where $d<a+1$ is a positive constant.  By Lemma $\ref{key lemma}$, we have

\begin{equation}\label{1}
\begin{aligned}
-\int_M & \nabla (\lambda ^b  f^a) \nabla f dv_g\\
&=\int_M \lambda^b f^a \Delta f dv_g\\
&\geq 2\int_M \lambda^b f^a |\nabla^{\perp}{\bf H}|^2 dv_g 
+2m\int_M \lambda^b f^{a+2} dv_g.
\end{aligned}
\end{equation}
On the other hand, we have

\begin{equation}\label{2}
\begin{aligned}
-\int_M & \nabla (\lambda ^b  f^a) \nabla f dv_g\\
&= -2 b \int_M \lambda^{b-1} \nabla \lambda f^a \langle \nabla^{\perp}{\bf H}, {\bf H} \rangle dv_g 
-4a \int_M \lambda^{b} f^{a-1} \langle \nabla^{\perp}{\bf H}, {\bf H} \rangle^2 dv_g \\
 & \leq  -2 b \int_M \lambda^{b-1} \nabla \lambda f^a \langle \nabla^{\perp}{\bf H}, {\bf H} \rangle dv_g \\
 & =  - 2 b \int_M \lambda^{b-1} \nabla \lambda f^a \langle \nabla^{\perp}{\bf H}, {\bf H}-\sqrt c X \rangle dv_g,
\end{aligned}
\end{equation}
where $X\in TM$ is a unit vector field.
From $(\ref{1})$ and $(\ref{2})$, we obtain

\begin{equation}\label{3}
\begin{aligned}
2\int_M & \lambda^b f^a |\nabla^{\perp}{\bf H}|^2 dv_g 
+2m\int_M \lambda^b f^{a+2} dv_g\\
 & \leq
- 2 b \int_M \lambda^{b-1} \nabla \lambda f^a \langle \nabla^{\perp}{\bf H}, {\bf H}-\sqrt c X \rangle dv_g\\
 & = -2  \int_M  \langle \lambda^{\frac{b}{2}} f^{\frac{a}{2}} \nabla^{\perp}{\bf H}, b \lambda^{\frac{b-2}{2}} f^{\frac{a}{2}}  \nabla \lambda ({\bf H}-\sqrt c X) \rangle dv_g \\
 & \leq
 \int_M \lambda^b f^a  |\nabla^{\perp}{\bf H}|^2 dv_g
 + \int_M b^2 \lambda^{b-2} f^a |\nabla \lambda|^2 f dv_g\\
 & =
 \int_M \lambda^b f^a  |\nabla^{\perp}{\bf H}|^2 dv_g
 + b^2 \int_M  \lambda^{b-2} f^{a+1} |\nabla \lambda|^2 dv_g.
\end{aligned}
\end{equation}

By using Young's inequality, we have

\begin{equation}\label{4}
\begin{aligned}
b^2 \int_M & \lambda^{b-2} f^{a+1} |\nabla \lambda|^2 dv_g\\
& = b^2 \int_M \lambda^{c} f^{d} \lambda^{b-2-c} f^{a+1-d} |\nabla \lambda|^2 dv_g \\
& \leq \int_M \lambda ^{b} f^{a+2} dv_g\\
& \ \ \ +C(a,d) \int_M \lambda ^{(a+1-d)\frac{a+2}{a+2-d}} f^{(a+1-d)\frac{a+2}{a+2-d}} |\nabla\lambda|^{2\frac{a+2}{a+2-d}} dv_g,
\end{aligned}
\end{equation}
where $c=\frac{d(a+3-d)}{a+2-d}$ and $C(a,d)$ is a constant depending only on $a$ and $d$.

Combining $(\ref{3})$ and $(\ref{4})$, we obtain
\begin{equation}\label{5}
\begin{aligned}
\int_M & \lambda^b f^a |\nabla^{\perp}{\bf H}|^2 dv_g 
+(2m-1)\int_M \lambda^b f^{a+2} dv_g\\
& \leq 
C(a,d) \int_M \lambda ^{(a+1-d)\frac{a+2}{a+2-d}} f^{(a+1-d)\frac{a+2}{a+2-d}} |\nabla\lambda|^{2\frac{a+2}{a+2-d}} dv_g \\
& \leq 
C(a,d) \int_M f^{(a+1-d)\frac{a+2}{a+2-d}} \left(\frac{1}{r}\right)^{2\frac{a+2}{a+2-d}} dv_g .
\end{aligned}
\end{equation}
Since $0<d<a+1$, we have $0<p:=(a+1-d)\frac{a+2}{a+2-d}<a+1$.
By the assumption $\int_M f^{ p }dv_g<\infty$ $(0<p<\infty)$, letting $r\nearrow \infty$ in $(\ref{5})$,
  the right hand side of $(\ref{5})$ goes to zero and the left hand side of $(\ref{5})$ goes to 
  $$\int_M  f^a |\nabla^{\perp}{\bf H}|^2 dv_g 
+(2m-1)\int_M  f^{a+2} dv_g,$$
   since $\lambda=1$ on $B_r(x_0)$.
 Thus, we have
  $f=|{\bf H}|^2-\left(c-\frac{\lambda}{2m}\right)=0,$ that is, the mean curvature is $\sqrt {c-\frac{\lambda}{2m}}$.
\end{proof}

\qquad\\

\section{Other results for BMO conjecture}\label{other}

The author introduced {\em polynomial growth bound of order $\alpha$ from below} as follows (cf. \cite{M2}).

\vspace{10pt}

\begin{defn}
For a complete Riemannian manifold $(N,h)$ and $\alpha \geq 0,$
if the sectional curvature $K^N$ of $N$ satisfies
$$K^N \geq-L(1+{\rm dist}_N(\cdot,q_0)^2)^{\frac{\alpha}{2}},~~{\rm for~some }~L>0~{\rm and}~q_0\in N,$$
then we say that $K^N$ has a {\em polynomial growth bound of order} $\alpha$ {\em from below.}
\end{defn}

\vspace{10pt}

An immersed submanifold $M$ in a Riemannian manifold $N$ is said to be {\em properly immersed} if the immersion is a proper map.
The author also showed as follows.

\vspace{10pt}

\begin{thm}[\cite{M2}]\label{Th of u}
Let $(M,g)$ be a properly immersed submanifold in a complete Riemannian manifold $(N,h)$ whose sectional curvature $K^N$ has a polynomial growth bound of order less than $2$ from below.
Let $u$ be a smooth non-negative function on $M$. 
Assume that there exists a positive constant $k > 0$ such that
\begin{equation}\label{key3 assumption}
\Delta u \geq k u^2\quad {\rm on}\ \ M.
\end{equation} 
If $|{\bf H}|\leq C(1+{\rm dist}_N(\cdot, q_0)^2)^{\frac{\beta}{2}}$ for some $C>0$, $0\leq\beta<1$ and $q_0\in N$, where ${\bf H}$ is the mean curvature vector field of $M$,
then $u = 0$ on $M$.
\end{thm}

\vspace{10pt}

By using Theorem $\ref{Th of u}$ and Lemma $\ref{key lemma}$, we obtain the following corollary.

\vspace{10pt}

\begin{cor}
Let $(N,h)$ be a complete Riemannian manifold with sectional curvature bounded from above by a non-negative constant $c$.  Assume that the curvature $K^N$ has a polynomial growth bound of order less than $2$ from below.
Let $\p:(M,g)\rightarrow(N,h)$ be a biminimal properly immersed submanifold with $\lambda\leq 2mc$ in $N$.
If $\sqrt {c-\frac{\lambda}{2m}} \leq|{\bf H}|\leq C(1+{\rm dist}_N(\cdot, q_0)^2)^{\frac{\beta}{2}}$ for some $C>0$, $0\leq\beta<1$ and $q_0\in N$, then the mean curvature is $\sqrt {c-\frac{\lambda}{2m}}$.
\end{cor}

\vspace{10pt}

Y. Luo introduced {\em at most polynomial volume growth} as follows (cf. \cite{yl1}).

\vspace{10pt}

\begin{defn}
Let $(M,g)$ be a complete Riemannian manifold and $x_0\in M$.
We say  $(M,g)$ is of {\em at most polynomial volume growth}, if there exists a non-negative integer $s$ such that 
$${\rm Vol}_g(B_{r}(x_0))\leq Cr^s,$$
where $C$ is a positive constant independent of $r$ and $B_{r}(x_0)$ is the ball centered at $x_0$ with radius $r$. 
\end{defn}

\vspace{10pt}

We obtain the following theorem.

\vspace{10pt}

\begin{thm}
Let $(N,h)$ be a Riemannian manifold with sectional curvature bounded from above by a non-negative constant $c$.
Let $\p:(M,g)\rightarrow (N,h)$ be a complete biminimal submanifold with $\lambda\leq2mc$ of at most polynomial volume growth in $N$.
 If $|{\bf H}|\geq \sqrt {c-\frac{\lambda}{2m}}$,
then the mean curvature is $\sqrt {c-\frac{\lambda}{2m}}$.
\end{thm}



\vspace{10pt}

\begin{proof}
Set $f=|{\bf H}|^2-\left(c-\frac{\lambda}{2m}\right)$. By the assumption $|{\bf H}|\geq \sqrt {c-\frac{\lambda}{2m}}$, $f$ is non-negative. 
For a fixed point $x_0\in M$, and for every $0<r<\infty,$
 we first take a cut off function $\lambda$ on $M$ satisfying that 
 \begin{equation}
\left\{
 \begin{aligned}
&0\leq\lambda(x)\leq1\ \ \ (x\in M),\\
&\lambda(x)=1\ \ \ \ \ \ \ \ \ (x\in B_r(x_0)),\\
&\lambda(x)=0\ \ \ \ \ \ \ \ \ (x\not\in B_{2r}(x_0)),\\
&|\nabla\lambda|\leq\frac{C}{r}\ \ \ \ \ \ \ (x\in M),\ \ \ \text{for some constant $C$ independent of $r$},
\end{aligned} 
\right.
\end{equation}
where $B_r(x_0)$ and  $B_{2r}(x_0)$ are the balls centered at a fixed point $x_0\in M$ with radius $r$ and $2r$ respectively (cf. \cite{N-U-2},\ \cite{M}).

Let $a$ be a positive constant to be determined later.  By Lemma $\ref{key lemma}$, we have

\begin{equation}\label{6}
\begin{aligned}
-\int_M & \nabla (\lambda ^{2a+4}  f^a) \nabla f dv_g\\
&=\int_M \lambda^{2a+4} f^a \Delta f dv_g\\
&\geq 2\int_M \lambda^{2a+4} f^a |\nabla^{\perp}{\bf H}|^2 dv_g 
+2m\int_M \lambda^{2a+4} f^{a+2} dv_g.
\end{aligned}
\end{equation}
On the other hand, we have

\begin{equation}\label{7}
\begin{aligned}
-\int_M & \nabla (\lambda ^{2a+4}  f^a) \nabla f dv_g\\
&= -2 ({2a+4}) \int_M \lambda^{2a+3} \nabla \lambda f^a \langle \nabla^{\perp}{\bf H}, {\bf H} \rangle dv_g \\
&\ \ \ -4a \int_M \lambda^{{2a+4}} f^{a-1} \langle \nabla^{\perp}{\bf H}, {\bf H} \rangle^2 dv_g \\
 & \leq  -2 ({2a+4}) \int_M \lambda^{2a+3} \nabla \lambda f^a \langle \nabla^{\perp}{\bf H}, {\bf H} \rangle dv_g \\
 & =  - 2 ({2a+4}) \int_M \lambda^{2a+3} \nabla \lambda f^a \langle \nabla^{\perp}{\bf H}, {\bf H}-\sqrt c X \rangle dv_g,
\end{aligned}
\end{equation}
where $X\in TM$ is a unit vector field.
From $(\ref{6})$ and $(\ref{7})$, we obtain

\begin{equation}\label{8}
\begin{aligned}
2\int_M & \lambda^{2a+4} f^a |\nabla^{\perp}{\bf H}|^2 dv_g 
+2m\int_M \lambda^{2a+4} f^{a+2} dv_g\\
 & \leq
- 2 ({2a+4}) \int_M \lambda^{2a+3} \nabla \lambda f^a \langle \nabla^{\perp}{\bf H}, {\bf H}-\sqrt c X \rangle dv_g\\
 & = -2  \int_M  \langle \lambda^{\frac{2a+4}{2}} f^{\frac{a}{2}} \nabla^{\perp}{\bf H}, ({2a+4}) \lambda^{\frac{2a+2}{2}} f^{\frac{a}{2}}  \nabla \lambda ({\bf H}-\sqrt c X) \rangle dv_g \\
 & \leq
 \int_M \lambda^{2a+4} f^a  |\nabla^{\perp}{\bf H}|^2 dv_g
 + \int_M ({2a+4})^2 \lambda^{2a+2} f^a |\nabla \lambda|^2 f dv_g\\
 & =
 \int_M \lambda^{2a+4} f^a  |\nabla^{\perp}{\bf H}|^2 dv_g
 + ({2a+4})^2 \int_M  \lambda^{2a+2} f^{a+1} |\nabla \lambda|^2 dv_g.
\end{aligned}
\end{equation}

By using Young's inequality, we have

\begin{equation}\label{9}
\begin{aligned}
({2a+4})^2 \int_M & \lambda^{2a+2} f^{a+1} |\nabla \lambda|^2 dv_g\\
& \leq \int_M \lambda ^{2a+4} f^{a+2} dv_g
 +C(a) \int_M  |\nabla\lambda|^{2(a+2)} dv_g,
\end{aligned}
\end{equation}
where $C(a)$ is a constant depending only on $a$.

Combining $(\ref{8})$,  $(\ref{9})$, we obtain
\begin{equation}\label{10}
\begin{aligned}
\int_M & \lambda^{2a+4} f^a |\nabla^{\perp}{\bf H}|^2 dv_g 
+(2m-1)\int_M \lambda^{2a+4} f^{a+2} dv_g\\
& \leq 
C(a) \int_M  |\nabla\lambda|^{2(a+2)} dv_g\\
& \leq 
C(a)  \left(\frac{1}{r}\right)^{2(a+2)} {\rm Vol}_g(B_{2r})\\
& \leq 
C(a)   r^{s-{2(a+2)}},
\end{aligned}
\end{equation}
where the last inequality follows from the assumption that there exists an integer $s\geq 0$ such that
$${\rm Vol}_g(B_{2r})\leq C r^s.$$

Choosing $a$ such that $a>{\rm max}\{0,\frac{s-4}{2}\}$.
Letting $r\nearrow \infty$ in $(\ref{10})$,
  the right hand side of $(\ref{10})$ goes to zero and the left hand side of $(\ref{10})$ goes to 
  $$\int_M  f^a |\nabla^{\perp}{\bf H}|^2 dv_g 
+(2m-1)\int_M  f^{a+2} dv_g,$$
   since $\lambda=1$ on $B_r(x_0)$.
 Thus, we have
  $f=|{\bf H}|^2-\left(c-\frac{\lambda}{2m}\right)=0,$ that is, the mean curvature is $\sqrt {c-\frac{\lambda}{2m}}$.
\end{proof}






\qquad\\
\qquad\\


\section{Biconservative hypersurfaces with $\lambda_i>-\frac{1}{2}mH$ in space forms}\label{sec bicon}

In this section, we consider the tangential part of a biharmonic hypersurface in a space form. It is called a {\em biconservative hypersurface} (cf.~\cite{CMOP2014}).

Let $\p:(M^m,g)\rightarrow (N^{m+1},h=\met)$ be an isometric immersion from an $m$-dimensional Riemannian manifold into an $m+1$-dimensional Riemannian manifold.
The Gauss formula is given by
\begin{equation}
\nabla^N_XY=\nabla _XY+B(X,Y),\ \ \ \ X,Y\in \frak{X}(M),
\end{equation}
where $B$ is the second fundamental form of $M$ in $N$.
The Weingarten formula is given by
\begin{equation}\label{2.Wformula}
\nabla^N_X \xi =-A_{\xi}X,\ \ \ X\in \frak{X}(M),
\end{equation}
where $A_{\xi}$ is the shape operator for the unit normal vector field $\xi$ on $M$.

For any $x \in M$, let $\{e_1, \cdots, e_m, e_{m+1}=\xi\}$ be an orthonormal basis of $N$ at $x$ such that $\{e_1, \cdots, e_m\}$ is an orthonormal basis of $T_xM$. 
Let $\lambda_1,\cdots, \lambda_m$ be the principal curvatures of a hypersurface $M$ at $x\in M$, that is, 
$A_{\xi}e_i=\lambda_i e_i$, for  $i=1,\cdots, m$.
Then, $B$ is decomposed as 
$$ 
B(X, Y) = B_{ m+1 }(X, Y)\xi,~~{\rm at}~x. 
$$ 
The squared norm of the second fundamental form is given by
$$|B|^2=|A|^2=\sum_{i=1}^m{\lambda_i}^2.$$
The mean curvature vector field ${\bf H}$ of $M$ at $x$ is also given by 
$$ 
{\bf H}(x) = \frac{1}{m} \sum_{i = 1}^m B(e_i, e_i) = H(x){\xi}=\frac{1}{m}\sum_{i=1}^m\lambda_i\xi,
$$ 
where $H(x) := \displaystyle\frac{1}{m} \displaystyle\sum_{i = 1}^m B_{ m+1 }(e_i, e_i)$
is the mean curvature.

By $(\ref{NS bih sub separate 2})$, a biconservative hypersurface in $N$ is defended as follows (cf.~\cite{CMOP2014}).

\vspace{10pt}

\begin{defn}[\cite{CMOP2014}]
An immersed hypersurface $\p:(M^m,g)\rightarrow (N^{m+1},h)$ is called {\em biconservative} if $\p:M^m\rightarrow N^{m+1}$ satisfies that
$$mH\nabla H +2A_{\xi}(\nabla H)-2H\left[{\rm Ric}^N(\xi)\right]^{T}=0,$$
where $H$ is the mean curvature, $\xi$ is the unit normal vector field of $M$ and ${\rm Ric^N}$ is the Ricci curvature of $N$.
 \end{defn}
 
 \vspace{10pt}
 
 From this, if $N$ is a space form, $\p:M\rightarrow N$ is biconservative if and only if
 \begin{align}\label{biconservative in spheres}
 mH\nabla H +2A_{\xi}(\nabla H)=0.
 \end{align}

We give an affirmative partial answer to BMO conjecture.
By using $(\ref{biconservative in spheres})$,
 we show as follows.

\vspace{10pt}
 
 \begin{prop}\label{prop bicon 1}
A biconservative hypersurface with $\lambda_i>-\frac{1}{2}mH$ $($for all $i=1,\cdots, m)$ in a space form has constant mean curvature.
 \end{prop} 

\vspace{10pt}
 
 \begin{proof}
 By $(\ref{biconservative in spheres})$, we have
\begin{align*}
 0=mH\nabla H +2A_{\xi}(\nabla H)
 &=mH\nabla H +2\sum_{i=1}^m\lambda_i(e_i H)e_i.
\end{align*}
Thus we have
$$(mH+2\lambda_j)e_jH=0, \ \ \ \ \ (\text{for all} ~~j=1, \cdots, m).$$
Therefore by the assumption, we have
$e_jH=0,$
at an arbitrary point $x\in M$ and for all $j=1, \cdots, m$. 
So we have that
$H$ is constant
at an arbitrary point $x\in M$.
Therefore $H$ is constant on $M$. 
 \end{proof}
 
 \vspace{10pt}
 
\begin{lem}\label{key lem of sphere}
A non-minimal biminimal hypersurface $M$ with constant mean curvature in a space form $N(c)$ with constant sectional curvature $c$ is a totally umbilical hypersurface with $|A|^2=cm-\lambda>0.$
\end{lem}

\vspace{10pt}

\begin{proof}
By $(\ref{biminimal eq})$, we have
$$\Delta H-H|A|^2+cmH-\lambda H=0.$$ 
Since $H$ is constant, 
\begin{align}\label{6.1}
H(|A|^2-cm+\lambda)=0.
\end{align}
Since $H\not=0,$ we have 
$$|A|^2=cm-\lambda.$$
Thus $cm-\lambda\geq 0$, but if $cm-\lambda=0$, we have $H=0$ which is a contradiction.
Therefore $cm-\lambda> 0$.
\end{proof}

From this and Proposition $\ref{prop bicon 1}$, we obtain the following result.

\vspace{10pt}

\begin{cor}\label{cor bicon 2}
If $M$ is a biharmonic hypersurface with $\lambda_i>-\frac{1}{2}mH$ $($for all $i=1,\cdots, m)$ in a space form $N(c)$ with constant sectional curvature $c$. 
Then we have $c>0$ and $M$ is a totally umbilical hypersurface with $|A|^2=cm.$
\end{cor}

\vspace{10pt}

\begin{proof}
From Proposition $\ref{prop bicon 1}$, we have that $H$ is constant.
By the assumption $\lambda_i>-\frac{1}{2}mH,$
we have 
 \begin{align}
H=\frac{1}{m} \sum_{i=1}^m\lambda_i>-\frac{1}{2}mH.
 \end{align}
Thus we have $H>0$.
By using Lemma $\ref{key lem of sphere}$, we obtain the result.
\end{proof} 
 
\vspace{10pt}
 
 \begin{rem}
If $M$ is a non-minimal biharmonic {\em convex} hypersurface, that is, a biharmonic hypersurface with the principal curvatures are non-negative, it satisfies that the condition $\lambda_i>-\frac{1}{2}mH$ $($see also \cite{yl2}$)$. 
 \end{rem}


\quad\\
\quad\\

\noindent 
{\bf Acknowledgements.} 
The author would like to express his gratitude to Professor Hajime Urakawa for many useful comments and valuable suggestions.

\bibliographystyle{amsbook}

\end{document}